\documentclass[10pt]{amsart}
\usepackage{amsmath,mathrsfs}
\usepackage{amssymb}
\usepackage{latexsym}
\usepackage{amsfonts}
\usepackage{hyperref}
\usepackage{graphicx}
\usepackage{psfrag}
\def\C{\mathbb C}

\def\R{{\mathbb R}}

\def\@tempb{saamsart}
\renewcommand{\theequation}{\thesection.\arabic{equation}}
\newtheorem{Pa}{Paper}[section]

\newtheorem{La}[Pa]{{\bf Lemma}}

\newtheorem{Pn}[Pa]{{\bf Proposition}}

\newtheorem{Ex}[Pa]{{\bf Example}}

\author[D. Alpay]{Daniel Alpay}
\address{(DA) Department of Mathematics
\newline
Ben Gurion University of the Negev \newline P.O.B. 653,
\newline
Be'er Sheva 84105, \newline ISRAEL} \email{dany@math.bgu.ac.il}
\author[I. Lewkowicz]{Izchak Lewkowicz}
\address{(IL) Department of Electrical Engineering
\newline
Ben Gurion University of the Negev \newline P.O.B. 653,
\newline
Be'er Sheva 84105, \newline ISRAEL }
\email{izchak@ee.bgu.ac.il}
\thanks{
This research was supported in part by the Bi-national Science
Foundation grant 2010117.}
\thanks{
D. Alpay thanks the
Earl Katz family for endowing the chair
which supported his research.}
\date{}
\title[Interpolating by polynomials with symmetries]
{interpolation by polynomials with symmetries
on the imaginary axis}

\begin{document}
\bibliographystyle{plain}
\begin{abstract}
We here specialize the standard matrix-valued polynomial
interpolation to the case where on the imaginary axis
the interpolating polynomials admit various symmetries:
Positive semidefinite, Skew-Hermitian, $J$-Hermitian,
Hamiltonian and others.

The procedure is comprized of three stages, illustrated
through the case where on $i\R$ the interpolating polynomials
are to be positive semidefinite. We first, on the expense
of doubling the degree, obtain a minimal degree interpolating
polynomial $P(s)$ which on $i\R$ is Hermitian. Then we find all
polynomials $\Psi(s)$, vanishing at the interpolation points
which are positive semidefinite on $i\R$. Finally, using the fact that
the set of positive semidefinite matrices is a convex subcone of
Hermitian matrices, one can compute the minimal scalar
$\hat{\beta}\geq 0$ so that $P(s)+\beta\Psi(s)$ satisfies all
interpolation constraints for all $\beta\geq\hat{\beta}$.

This approach is then adapted to cases when the family of
interpolating polynomials is not convex. 
Whenever convex, we
parameterize all minimal degree interpolating polynomials.
\end{abstract}
\maketitle
\renewcommand{\theequation}
{\thesection.\arabic{equation}}

\section{introduction}\label{sec:introduction}
\setcounter{equation}{0}

Probably the simplest version of interpolation problem is as
follows. Given a family of functions $\mathcal{F}$, nodes
$x_1~,~\ldots~,~x_p$ and image points $Y_1~,~\ldots~,~Y_p$, search
for $F\in\mathcal{F}$ so that
\begin{equation}\label{eq:BasicInterp}
Y_j=F(x_j)\quad\quad\quad\quad j=1~,~\ldots~,~p.
\end{equation}
More specifically, find out whether such $F(s)$ exists and if yes,
search for all ``simple" interpolating functions in $\mathcal{F}$.
In the context of rational functions, ``simple" means low degree,
which here takes the form of the McMillan degree, see e.g.
\cite{AV}, \cite{BGR}, \cite{BK}, \cite{GL}.
\vskip 0.2cm

There is a vast literature on this classical problem. For a
comprehensive study see e.g. \cite{BGR}. Additional relevant
literature is presented in Section \ref{sec:InterpolationBackground}.
\vskip 0.2cm

Using the framework in \eqref{eq:BasicInterp}, in this work, the
nodes $x_j$ are in $\C$, the image points $Y_j$ are $m\times m$
matrices and $\mathcal{F}$ are polynomials $F$ of a complex
variable $s$, i.e.
\begin{equation}\label{eq:Polynom}
F(s)=\sum\limits_{k=0}^qC_ks^k\quad\quad\quad\quad
C_k\in\C^{m\times m}.
\end{equation}
Recall, that the McMillan degree of matrix-valued polynomials
is well defined. In \cite[Corollary 2.1.1]{BK} it is shown that 
for $F(s)$ in \eqref{eq:Polynom} it is equal to the rank of the
block-triangular, block-Toeplitz matrix
\begin{equation}\label{eq:Deg}
\left(\begin{smallmatrix}
C_q&C_{q-1}&C_{q-2}&\ldots&C_1   &C_o   \\
0      &C_q&C_{q-1}&\ldots&C_2   &C_1   \\
0      & 0     &C_q&\ldots&C_3   &C_2   \\
\vdots &\vdots &\vdots &\vdots&\vdots&\vdots\\
0      & 0     & 0     &\ldots&  0   &C_q
\end{smallmatrix}\right).
\end{equation}
In particular, if $C_q$ is nonsingular, then the McMillan degree of
$F(s)$ is equal to $m(q+1)$.
\vskip 0.2cm

In this work, the polynomials in $\mathcal{F}$ are restricted to
have various symmetries on the imaginary axis, described in the
sequel. The importance of polynomial matrix interpolation with
symmetries, was raised in \cite[Subsection 2.1.48]{Hi}. In the
framework of~ {\em real}~ variable, this problem has already been
treated in \cite{DMPP}. That work differs from ours in many ways.
\vskip 0.2cm

Through Examples \ref{Ex:Motivation1}, \ref{Ex:Motivation2},
\ref{Ex:Motivation3}, \ref{Ex:Motivation4} and \ref{Ex:Psia} part
A, we illustrate the fact that even in the scalar case, the question
addressed here is not trivial. We start with the following.

\begin{Ex}\label{Ex:Motivation1}
{\rm

As a prototype example we shall seek a ~{\em minimal degree}~
interpolating function $F(s)$ so that
\[
\begin{smallmatrix}
F(s)&\left(\begin{smallmatrix}1\\2\\3\end{smallmatrix}\right)&
\rightarrow&
\left(\begin{smallmatrix}18\\75\\50\end{smallmatrix}\right)
\end{smallmatrix}.
\]
Here and in Examples \ref{Ex:Motivation2}, \ref{Ex:Motivation3}
and \ref{Ex:Motivation4}, we shall consider various families
of functions $\mathcal{F}$ with the same points.
\vskip 0.2cm

Taking the family $\mathcal{F}$ to be ~{\em unstructured}~
polynomials a straightforward computation yields,
%
\[
F_1(s)=-45s^2+180s-121.
\]
}
\qed
\end{Ex}

As mentioned, we here focus on ~{\em structured}~ polynomials. To
formally set-up the problem addressed here, we need some background.

\subsection{Functions with symmetry on
the imaginary axis}\label{sec:Functions}

Let $F$ be $m\times m$-valued rational function in the sense
that
\[
F:\C~\rightarrow~\C^{m\times m}.
\]
In the sequel we shall use the following notation,
\[
F^{\#}(s):=F^*(-s^*).
\]
We shall denote by $\C_+$ $(\overline{\C_+})$ the open (closed)
right half plane and by
$\overline{\mathbb P}_m$ $({\mathbb P}_m)$ the sets
of $m\times m$ positive semidefinite (definite) matrices. Whenever
clear from the context, the subscript will be omitted and we
shall simply write $\overline{\mathbb P}$ or ${\mathbb P}$.
\vskip 0.2cm

We call functions $F(s)$ ~{\em Positive},~ denoted by
$\mathcal{P}$, if they are analytic in $\C_+$ and
\begin{equation}\label{eq:pos}
\left(F(s)+F^*(s)\right)\in\overline{\mathbb P}
\quad\quad\quad s\in\C_+~.
\end{equation}
These functions has played an important in the theory of
electrical networks from around 1930, see e.g. \cite{AV} and
\cite{Be}. They also serve as the corner stone of the theory
of linear dissipative systems (a.k.a absolutely stable), see
e.g. \cite[Theorem 2.7.1]{AV}, \cite[3.18]{Be}.
\vskip 0.2cm

One can relax the condition and call functions $F(s)$~
{\em Generalized Positive}, denoted by $\mathcal{GP}$, if
\[
\left(F(s)+F^*(s)\right)\in\overline{\mathbb P}
\quad\quad\quad\quad
{\rm almost~for~all}\quad s\in{i}\R.
\]
For an early study of rational $\mathcal{GP}$ functions see
\cite{AM} and for recent references see e.g. \cite{AL1},
\cite{AL2} and \cite{AL3}. These functions were studied in
different frameworks, for example when the upper half plane
replaces $\C_+$, see e.g.  \cite{DHdS}, \cite{DLLS1},
\cite{L1} or when the unit disk replaces $\C_+$, either in
the
domain see e.g. \cite{DGK} or the image, see e.g. \cite{GR}.
\vskip 0.2cm

Abusing terminology of real scalar functions, we call a
function $F(s)$~ {\em Odd}~ if
\begin{equation}\label{eq:Odd}
F^{\#}(s)=-F(s).
\end{equation}
This implies that on the imaginary axis $F(s)$ is
skew-Hermitian, i.e.
\[
\left(F(s)_{|_{s\in{i}\R}}\right)^*=-F(s)_{|_{s\in{i}\R}}~.
\]
We shall denote by ${\mathcal Odd}$ the set of odd
functions. Note that ${\mathcal Odd}\subset\mathcal{GP}$,
for details, see \cite[Proposition 4.2]{AL3}.
\vskip 0.2cm

In a similar way we call $F(s)$ ~{\em Even}~ if
\begin{equation}\label{eq:even}
F^{\#}(s)=F(s).
\end{equation}
This implies that
on the imaginary axis $F(s)$ is Hermitian, i.e.
\[
\left(F(s)_{|_{s\in{i}\R}}\right)^*=F(s)_{|_{s\in{i}\R}}~.
\]
The set of even functions, denote by ${\mathcal Even}$, was
studied in \cite[Section 5]{AL3}.
\vskip 0.2cm

Of particular interest is the class of ~{\em Generalized Positive
Even}, denoted by $\mathcal{GPE}$,
\begin{equation}\label{eq:DefGPE}
\begin{smallmatrix}
F^{\#}(s)=F(s)&~&~&~\\~\\
F(s)\in\overline{\mathbb P}&~&~&
{\rm almost~for~all}\quad s\in{i}\R.
\end{smallmatrix}
\end{equation}
For details see \cite[Section 5]{AL3}. Recall that
$F\in\mathcal{GPE}$ if and only if there exist $G(s)$ so that
\begin{equation}\label{SpectFact}
F(s)=G(s)G^{\#}(s).
\end{equation}
If $F(s)$ is analytic on the imaginary axis \eqref{SpectFact}
is called ~{\em spectral factorization}\begin{footnote}{$G(s)$
is analytic in $\overline{\C_+}$
sometimes also $G^{-1}(s)$ is analytic there.}\end{footnote},
see e.g. \cite[Section 5.2]{AV}, \cite[Chapter 9]{BGKR},
and \cite[Section 19.3]{LR} 
Else, \eqref{SpectFact} is a ~{\em pseudo spectral factorization},
see e.g. \cite[Chapter 10]{BGKR},
\vskip 0.2cm

For future reference we
recall that a convex cone which in addition is closed under
inversion is called a Convex Invertible Cone, {\bf cic}~ in
short\begin{footnote}{Strictly speaking, this means that
whenever the inverse exists, it also belongs to the set, e.g.
the set of positive semidefinite matrices is a {\bf cic}. In
contrast, the open upper half of $\C$ is not.}
\end{footnote}, see e.g. \cite{CL1}, \cite{CL2} and \cite{CL4}.
\vskip 0.2cm

It is easy to see that the sets ${\mathcal Odd}$ and
${\mathcal Even}$ are closed under positive scaling, summation
and inversion, i.e. {\bf cic}s.
\vskip 0.2cm

The set ${\mathcal P}$ is a sub{\bf cic} of $\mathcal{GP}$ functions.
More precisely, ${\mathcal P}$ is a maximal {\bf cic}~ of functions
which are analytic in $~\C_+$, see e.g. \cite[Proposition 4.1.1]{CL4}.
The {\bf cic} structure of $\mathcal{GP}$ functions was studied
in \cite{AL3}.
\vskip 0.2cm

Recall that the intersection of {\bf cic}s is a {\bf cic}, e.g.
\cite[Proposition 2.2]{CL1}. In particular $\mathcal{GPE}$ is
the intersection between the {\bf cic}s of $\mathcal{GP}$ and
${\mathcal Even}$. We now describe an intermediate set between
the ${\mathcal Even}$ its sub{\bf cic} $\mathcal{GPE}$.
Recall that ${\mathcal Even}$, see \eqref{eq:even}, means
$F_{|_{s\in{i}\R}}$ is
Hermitian and $\mathcal{GPE}$ means $F_{|_{s\in{i}\R}}$ is
positive semidefinite, see \eqref{eq:DefGPE}.
\vskip 0.2cm

Next consider the set where for all $\omega\in\R$, $F(i\omega)$
is Hermitian and there is no eigenvalues crossing from $\C_-$
to $\C_+$ (or vice versa). Roughly, along the imaginary axis
the inertia is (almost) fixed. More precisely, we call a
$m\times m$-valued rational function $F(s)$~
{\em $\nu$-Generalized Positive Even},~ denoted by
$\mathcal{{\nu}GPE}$ if it admits a factorization
\[
F(s)=G(s){\rm diag}\{-I_{\nu}~,~I_{m-\nu}\}G^{\#}(s)
\quad\quad\quad\quad \nu\in[0,~l],
\]
for some $G(s)$. In particular, for $\nu=0$ one returns to
the $\mathcal{GPE}$ case in \eqref{eq:DefGPE} and if
$\nu=m$, then $-F\in\mathcal{GPE}$. Thus, to avoid
triviality,
in the sequel we shall focus on $\mathcal{{\nu}GPE}$
functions admitting factorization
\begin{equation}\label{JSpectFact}
F(s)=G(s){\rm diag}\{-I_{\nu}~,~I_{m-\nu}\}G^{\#}(s)
\quad\quad\quad\quad \nu\in[1,~m-1].
\end{equation}
\vskip 0.2cm

If $F(s)$ is analytic on the imaginary axis \eqref{JSpectFact}
is called ~{\em $J$-spectral factorization}.
Else, this is $J$-pseudo spectral factorization, see e.g.
\cite[Part VII]{BGKR}.
\vskip 0.2cm

It should be pointed out that the technique we use in the sequel
does not require factorization of $\mathcal{{\nu}GPE}$ (or
$\mathcal{GPE}$) functions.
\vskip 0.2cm

{\em
The aim of this work is to offer a way for solving the interpolation
problem in \eqref{eq:BasicInterp} where the family $\mathcal{F}$
is comprized of matrix-valued polynomials within: ${\mathcal Odd}$,
${\mathcal Even}$, $\mathcal{GPE}$ and $\mathcal{{\nu}GPE}$.
Moreover, whenever $\mathcal{F}$ is convex, all minimal degree
interpolating polynomials are given.
}
\vskip 0.2cm

We now set up the main idea of this work in terms of
\eqref{eq:BasicInterp}, through the framework of
$\mathcal{F}=\mathcal{GPE}$ polynomials.
Other cases will turn to be (not necessarily small) variations
on the same theme.

\begin{itemize}
\item[(i)~~~]{}In classical unstructured polynomial interpolation
one obtains $P(s)$ ($P~:~\C~\rightarrow~\C^{m\times m}$) 
\begin{equation}\label{P}
P(s)=\sum\limits_{k=0}^{q-1}C_ks^{k}\quad\quad\quad
C_k\in\C^{m\times m},
\end{equation}
an unstructured interpolating polynomial of minimal degree,
namely it is at most $mq$ with $p-1\geq q$ (recall $p$ is the
number of interpolating nodes).

\noindent
On the expense of doubling its degree (i.e.
$2m(n-1)\geq{\rm deg}(P)$), one can introduce partial structure to
the interpolating polynomial $P(s)$ in \eqref{P}, see Proposition
\ref{Pr:Lagrange}. In particular, the resulting polynomial
$P(s)$ may be in ${\mathcal Even}$, see Subsection
\ref{sec:EvenLagrange}, or in ${\mathcal Odd}$, see Subsection
\ref{sec:Odd}.

\noindent
If the obtained $P(s)$ is already in $\mathcal{GPE}$, we are
done. Assume that this is not the case.

\item[(ii)~~]{}It is easy to obtain all minimal degree $\mathcal{GPE}$
polynomials vanishing at $x_1~,~\ldots~,~x_n$, denoted by $\Psi(s)$,
$(2mn\geq{\rm deg}(\Psi)$), see Section \ref{sec:NeutPoly} and
Proposition \ref{Pn:GpeNeut}.

\item[(iii)~]{}Using the above $P(s)$ and $\Psi(s)$, for all
$\beta\in\C$,
\begin{equation}\label{eq:InterpPolynom}
F(s)=P(s)+\beta\Psi(s),
\end{equation}
is an interpolating polynomial. In fact, $F(s)$ is in
${\mathcal Even}$ for all $\beta\in\R$.

\item[(iv)~]{}The family of $\mathcal{GPE}$ polynomials is a
convex subcone of ${\mathcal Even}$ and by construction
\mbox{${\rm deg}(\Psi)>{\rm deg}(P)$}. Thus, one can find
$\hat{\beta}$ so that for all $\beta\geq\hat{\beta}$ the
interpolating polynomial $F(s)$ in \eqref{eq:InterpPolynom}
is in $\mathcal{GPE}$. See Proposition \ref{Pn:HatBeta}.
\end{itemize}
\vskip 0.2cm

As mentioned, interpolation with $\mathcal{GPE}$ polynomials
serves as a prototype of our technique. In Section
\ref{sec:SpecialCases} the same idea is extended to polynomials
with various symmetries on the imaginary axis. In Section
\ref{sec:refine} refine the procedure by exploiting
possible special structure of $P(s)$ obtained in stage (i) of
the ``recipe" to enlarge the family of polynomials $\Psi(s)$
obtained in stage (ii) of the ``recipe" and to reduce the
minimal degree of the interpolating polynomial $F(s)$.
In Section \ref{sec:J-GPE} we modify the recipe to allow
interpolation with non-convex set of polynomials, e.g.
$\mathcal{{\nu}GPE}$ in \eqref{JSpectFact}. Finally, in Section
\ref{sec:FutureResearch} a sample of future research problems.

\section{motivation and background}
\label{sec:InterpolationBackground}
\setcounter{equation}{0}

In this section we review various aspects of the interpolation
problem at hand. For completeness we start by presenting 
three of the popular variants of the interpolation problem
not addressed in this work.
\vskip 0.2cm

First, the ~{\em tangential interpolation}, i.e.
given 
$x_1,~\ldots~,~x_p\in\C$,~
$v_1,~\ldots~,~v_p\in\C^m$
and
$Y_1,~\ldots~,~Y_p\in\C^l$
find $F(s)$ in $\mathcal{F}$ $(F~:~\C~\rightarrow~\C^{l\times m})$
so that
\begin{equation}\label{eq:RightTangential}
Y_j=F(x_j)v_j \quad\quad\quad\quad j=1,~\ldots~,~p.
\end{equation}
Strictly speaking, this is the ~{\em right}~ tangential
interpolation problem. The ~{\em left}~ tangential interpolation
problem is:~ For given $x_1,~\ldots~,~x_p\in\C$,~
$u_1,~\ldots~,~u_n\in\C^l$
and $Z_1,~\ldots~,~Z_p\in\C^m$ finding $F(s)$
$(F~:~\C~\rightarrow~\C^{l\times m})$
so that
\begin{equation}\label{eq:LeftTangential}
Z_j^*=u_j^*F(x_j)\quad\quad\quad\quad j=1,~\ldots~,~p.
\end{equation}
Combining both, one obtains the ~{\em bi-tangential}~ interpolation
problem. For sample references see e.g. 
\cite{BK}, \cite{Fu3} for the polynomial case, 
\cite{ABKW} for the unstructured rational case,
and \cite{ABL1}, \cite{ABL2}, 
for the structured rational case.
\vskip 0.2cm

As a second popular variant we mention that in control theory
there has been an interest in problems of the following form.
Given $A(s), B(s)$ $m\times m$-valued polynomials where in
addition $B\in{\mathcal Even}$ one seeks all $m\times m$-valued
polynomials $X(s)$ satisfying,
\[
X^{\#}(s)A(s)+A^{\#}(s)X(s)=2B(s).
\]
See e.g. \cite{HS} and the survey in \cite{CZ}
\vskip 0.2cm

To present a third popular variant, we start by resorting to the
notion of ~{\em reverse polynomial}. Recall that $\hat{F}(s)$ is
said to be the reverse of a polynomial
\mbox{$F(s)=\sum\limits_{j=0}^qs^jC_j$} in \eqref{eq:Polynom} if
\begin{equation}\label{eq:reverse}
\hat{F}(s):=s^qF(s^{-1})=\sum\limits_{j=0}^qs^jC_{q-j},
\end{equation}
see e.g. \cite[Section 2]{ACL}, \cite[Eq. (8.17)]{Fu2}.
\vskip 0.2cm

Next, we recall in the notion of ~{\em linearization}.~ For
constant $mq\times mq$ matrices $A, B$ we say that the pencil
$sA-B$ is a {\em linearlization}~ of the polynomial $F(s)$
in \eqref{eq:Polynom} if there exist $mq\times mq$-valued
unimodular polynomials $E(s), D(s), G(s), H(s)$ so that,
\[
\begin{matrix}
{\rm diag}\{F(s),&I_{m(q-1)}\}&=&E(s)(sA-B)D(s)\\
{\rm diag}\{\hat{F}(s),&I_{m(q-1)}\}&=&G(s)(A-sB)H(s).
\end{matrix}
\]
For details see e.g. \cite{ACL}.
\vskip 0.2cm

Roughly speaking, out of the many variants of interpolations,
two of the better studied frameworks are (i) unstructured
polynomials of the Lagrange type see e.g.
\cite{ACL}, \cite{BT}, \cite{La}, \cite{Mej}
and (ii) structured rational
functions, of the Nevanlinna-Pick type, see e.g.
\cite[Section 18]{BGR}.
\vskip 0.2cm

There are fundamental differences between these two problems:
The Nevanlinna-Pick type problems is far
more involved, but the underlying structure allows for the
use of powerful tools.
\vskip 0.2cm

Addressing interpolation through structured polynomials, as
we do here, turn most of the classical arsenal, like linear
fractional transformation, irrelevant.
\vskip 0.2cm

Nevanlinna-Pick interpolation of positive functions $\mathcal{P}$
has been well studied, see e.g. \cite[Chapter 18]{BGR}.
It was extended, not in the framework of $\mathcal{GP}$ but of
(i) generalized Schur functions (contractive on the unit circle),
in numerous works, see e.g. \cite{AADLW}, \cite{Ball}, \cite{BH},
\cite{bol_oam2}, and \cite{DD} and (ii) generalized Nevanlinna
functions (mapping the real axis to the upper half plane)
\cite{ABD}, \cite{AmDe}, \cite{ADLS}, \cite{bol_oam},
\cite[Section 3]{DZ} and \cite{ADLRS}.
\vskip 0.2cm

It should be pointed out that this extension of Nevanlinna-Pick
interpolation from addressing $\mathcal{P}$ to $\mathcal{GP}$
is computationally involved. Moreover, the
existing parameterization of all $\mathcal{GP}$ interpolating
functions (after being translated from the generalized Schur
framework) neither single out $\mathcal{GPE}$ functions nor
polynomials. This is illustrated next.

\begin{Ex}\label{Ex:Motivation2}
{\rm
Consider the interpolation points in Example \ref{Ex:Motivation1}.
\vskip 0.2cm

I. Assume now that the family $\mathcal{F}$ is the set of
rational $\mathcal{GP}$ functions. Then, minimal degree
interpolating function is
\[
F_2(s)=\frac{150(574-451s)}{41(66-41s)}~.
\]
Clearly, $F_2(s)$ is neither in ${\mathcal Even}$ nor a polynomial.
\vskip 0.2cm

II. Consider the interpolating polynomial in Example
\ref{Ex:Motivation1}. Then, $F_1(s)$ is neither in ${\mathcal Even}$
nor in $\mathcal{GP}$.
}
\qed
\end{Ex}

To summarize, from a practical point of view, the existing
interpolation scheme for $\mathcal{GP}$ functions, is not very
helpful for interpolation by $\mathcal{GPE}$ polynomials.
\vskip 0.2cm

{\bf The Lagrange approach}

A classical approach, see e.g. \cite[Section 2.10]{Fu2}, common
to some interpolation problems of the
form \eqref{eq:BasicInterp} is here illustrated
\vskip 0.2cm

Take,
\[
F(s)=\sum\limits_{j=1}^p\tilde{F}_j(s)\quad\quad\quad\quad
\quad\quad\quad\quad
\tilde{F}_j
\in\mathcal{F},\quad j=1,~\ldots~,~p
\]
so that
\[
\begin{smallmatrix}~&~&x_1   &~&x_2   &~&\ldots&~&x_p   \\
\tilde{F}_1(s)      &~&Y_1   &~&0     &~&\ldots&~&0     \\
\tilde{F}_2(s)      &~&0     &~&Y_2   &~&\ldots&~&0     \\
\vdots              &~&\vdots&~&\vdots&~&\cdots&~&\vdots\\ 
\tilde{F}_p(s)      &~&  0   &~& 0    &~&\ldots&~&Y_p.
\end{smallmatrix}
\]
In the framework of scalar (non-structured) polynomial interpolation,
this approach probably preceded Lagrange \cite{La}
(in \cite{Mej} it is attributed to \cite{Wa}).
\vskip 0.2cm

Computationally motivated, an interesting choice of $\tilde{F}_j(s)$
for unstructured scalar polynomials was presented in \cite[Eq. (4.2)]{BT}.
\vskip 0.2cm

The above straightforward approach to interpolation problems has
some limitations:

(i) It first assumes that it is easy
to construct the elements $\tilde{F}_j
$
in $\mathcal{F}$ (for example if $\mathcal{F}$ is the set of
$\mathcal{GP}$ polynomials, this is not easy, see e.g.
\cite[Example 5.3b]{AL3}).

(ii) It assumes that $\mathcal{F}$, the family of interpolating
functions, is convex (in Section \ref{sec:J-GPE} we address a
non-convex family $\mathcal{F}$).

(iii) In addition if for example $\mathcal{F}$ is a convex set of
rational functions this scheme may yield high degree
interpolating functions.
\vskip 0.2cm

We now further scrutinize this interpolation scheme. To this
end we here employ it to an example\begin{footnote}{We have
already used it in the framework of $\mathcal{GPE}$
polynomials in \cite[Example 5.3a]{AL3}.}\end{footnote}.
\vskip 0.2cm

\begin{Ex}\label{Ex:Motivation3}
{\rm

We here show that adapting the Lagrange approach to the case
where $\mathcal{F}=\mathcal{GPE}$ polynomials, enables us to
construct ~{\em some}~ of the~ {\em minimal degree}~ interpolating
polynomials in the problem addressed in Examples
\ref{Ex:Motivation1} and \ref{Ex:Motivation2}.

Indeed, take
\[
F_3(s)
=\tilde{F}_1(s)+\tilde{F}_2(s)+\tilde{F}_3(s)
\]
with
\[
\begin{smallmatrix}
\tilde{F}_1(s)&=&
(4-s^2)(9-s^2)(\tilde{\alpha}(1-s^2)+\frac{3}{4})&~&
\tilde{\alpha}\geq 0,\\
\tilde{F}_2(s)&=&
(1-s^2)(9-s^2)(\tilde{\beta}(4-s^2)-\frac{5}{4}s^2)
&~&\tilde{\beta}\geq 0,\\
\tilde{F}_3(s)&=&
(1-s^2)(4-s^2)(\tilde{\gamma}(9-s^2)+\frac{5}{4})
&~&\tilde{\gamma}\geq 0.\end{smallmatrix}
\]
This can be aggregated to a single parameter
\begin{equation}\label{eq:ExPsib}
F_3(s)=-3s^4+34s^2-13+
\beta
(1-s^2)(4-s^2)(9-s^2)\quad\quad\quad
\beta\geq\frac{5}{4}~.
\end{equation}
}
\qed
\end{Ex}

%
In Example 
%
\ref{Ex:Motivation4} 
below we show that
in contrast to to \eqref{eq:ExPsib}, there is a $\mathcal{GPE}$ 
interpolating polynomial for all $\beta\geq\frac{1}{2}~$. Thus,
adapting
the Lagrange approach to $\mathcal{GPE}$ polynomials is
appealing due to its simplicity. However, even in the scalar
case it turns out to provide conservative results. See also
Part A of Example \ref{Ex:Psia}.
\vskip 0.2cm

We conclude this section by further examining
the
structure of some the families $\mathcal{F}$ involved.
Let
$F(s)$ be a $m\times m$-valued polynomial as in \eqref{eq:Polynom}.
It is straightforward to verify that
\[
\begin{smallmatrix}
\mathcal{F}    &~&~&C_{2k}             &~&C_{2k+1}            \\~\\
{\mathcal Even}&~&~&{\rm Hermitian}    &~&{\rm skew-Hermitian}\\~\\
{\mathcal Odd}&~&~&{\rm skew-Hermitian}&~&{\rm Hermitian}.
\end{smallmatrix}
\]
However that there is no explicit way to characterize the sets
$\mathcal{GP}$, $\mathcal{GPE}$ $\mathcal{{\nu}GPE}$ only
through the structure of their coefficients. For example,
scalar, second degree, $\mathcal{GPE}$ polynomials are given by
\[
a+b(r+is)^2\quad\quad{\rm with}\quad\quad a\geq 0\quad
b>0\quad r\in\R
\]
(it is a strict subset of all $F(s)$ in \eqref{eq:Polynom} with
$q=2$ where $C_o\geq 0$, $C_1\in{i}\R$ and $0>C_2$). Hence, it is
conceivable to presume that to solve the problem in
\eqref{eq:BasicInterp} one needs to go beyond a simple modification
of the classical polynomial
interpolation.

\section{partly structured polynomial interpolation}
\label{sec:Lagrange}
\setcounter{equation}{0}

The classical $m\times m$-valued unstructured polynomial
interpolation can be formulated as follows. Given
\mbox{$x_1,~\ldots~,~x_p\in\C$}
(distinct) and \mbox{$Y_1~,~\ldots~,~Y_p\in\C^{m\times m}$}
find $P(s)$, a minimal degree polynomial
($P~:~\C~\rightarrow~\C^{m\times m}$) see \eqref{P}, 
so that
\[
Y_j=P(x_j)\quad\quad\quad\quad j=1,~\ldots~,~p.
\]
It is
known that the problem is solvable and in \eqref{P}, $p\ge q$.
\vskip 0.2cm

We here adapt\begin{footnote}{Originally, it appeared in other
frameworks: Carath\'{e}odory functions in \cite{ABL1} and the unit
disk $H_2$ functions in \cite{ABL2}.}\end{footnote} an idea
from \cite{ABL1}, \cite[Section 2]{ABL2} enabling us, by
roughly doubling $q$ in \eqref{P}, to impose structure on the
matrical coefficients $C_k$. The problem can be formulated as
follows.

\vskip 0.2cm
Given: $x_1,~\ldots~,~x_p\in\C$ and
\mbox{$A, B, Y_1~,~\ldots~,~Y_p\in\C^{m\times m}$}
find $P(s)$ a low degree $m\times m$-valued polynomial
($P~:~\C~\rightarrow~\C^{m\times m}$) so that
\begin{equation}\label{InterpP}
\begin{smallmatrix}
P^{\#}(s)&=&AP(s)B&~&~&A, B\in\C^{m\times m}\quad
{\rm non-singular}\\~\\
P(x_j)&=&Y_j&~&~&j=1,~\ldots~,~p.
\end{smallmatrix}
\end{equation}
To guarantee feasibility of the problem one needs to assume that
the data satisfies
\begin{equation}\label{eq:feasible}
\begin{smallmatrix}
x_j-x_k=0&~&\Longrightarrow&~&Y_j&=&Y_k\\~\\
x_j+x_k^*=0&~&\Longrightarrow&~&Y_j&=&(AY_kB)^*
\end{smallmatrix}\quad\quad\quad p\geq k>j\geq 1.
\end{equation}
Next, we shall call the interpolation data set~ {\em reduced}
if out of feasible points $x_1~,~\ldots~,~x_p$ we extract a
maximal subset $x_1~,~\ldots~,~x_n$, i.e. $p\geq n$ (and
the corresponding $Y_1,~\ldots~,~Y_n$) so that
\begin{equation}\label{eq:RedData}
\begin{smallmatrix}
x_j-x_k\not=0\\~\\
x_j+x_k^*\not=0
\end{smallmatrix}\quad\quad\quad n\geq k>j\geq 1.
\end{equation}
Note that the choice of $x_1~,~\ldots~,~x_n$ out of the $p$
given nodes, is not unique ($n$ is unique) but it will not
affect the proposed procedure below\begin{footnote}{To summarize,
$m\times m$ is the dimension of $F(s)$, the number of given
interpolation points is $p$ and it is then reduced to $n$.
Finally $q$ is the number of matrical coefficients in
\eqref{P} and \eqref{eq:InterpPolynom}.}\end{footnote}.
\vskip 0.2cm

\begin{Pn}\label{Pr:Lagrange}
Given a feasible data set \eqref{eq:feasible}. There always
exists an interpolating polynomial $P(s)$ \eqref{P} satisfying
\eqref{InterpP} with
\[
2n\geq q,
\]
where $n$ is the dimension of the reduced data set
\eqref{eq:RedData}.

The coefficients $C_o,~\ldots~,~C_{2n-1}$ are obtained from the
following matrix equation
\begin{equation}\label{eq:Lagrange}
\mbox{\large\boldmath{X}}
\mbox{\large\boldmath{C}}=
\mbox{\large\boldmath{Y}}
\end{equation}
where the dimensions of both $\mbox{\large\boldmath{C}}$
and $\mbox{\large\boldmath{Y}}$ is $2nm\times m$,
\[
\mbox{\large\boldmath{C}}:=
\left(\begin{smallmatrix}C_o\\ \vdots\\ C_{2n-1}
\end{smallmatrix}\right)\quad\quad\quad
\mbox{\large\boldmath{Y}}
:=\left(\begin{smallmatrix}Y_1\\ \vdots\\ Y_n\\
(AY_1B)^*\\ \vdots\\ (AY_nB)^*\end{smallmatrix}\right)
\]
and $\mbox{\large\boldmath{X}}$ is
the $2nm\times 2nm$ block-Vandermonde matrix
\[
\mbox{\large\boldmath{X}}=
\left(
\begin{smallmatrix}
I_m&x_1I_m&x_1^2I_m&\ldots&x_1^{2n-1}I_m\\
I_m&x_2I_m&x_2^2I_m&\ldots&x_2^{2n-1}I_m\\
\vdots&\vdots&\vdots&\vdots&\vdots      \\
I_m&x_nI_m&x_n^2I_m&\ldots&x_n^{2n-1}I_m\\
I_m&-x_1^*I_m&(-x_1^*)^2I_m&\ldots&(-x_1^*)^{2n-1}I_m\\
I_m&-x_2^*I_m&(-x_2^*)^2I_m&\ldots&(-x_2^*)^{2n-1}I_m\\
\vdots&\vdots&\vdots&\vdots&\vdots      \\
I_m&-x_n^*I_m&(-x_n^*)^2I_m&\ldots&(-{x_n}^*)^{2n-1}I_m
\end{smallmatrix}\right).
\]
\end{Pn}

{\bf Proof}\quad
Note that $\mbox{\large\boldmath{X}}$ 
can be written as
\[
\mbox{\large\boldmath{X}}=\hat{X}\otimes{I_m},
\]
where $\hat{X}$ is an actual $2n\times 2n$ Vandermonde matrix (i.e.
substitute $m=1$ in $\mbox{\large\boldmath{X}}$
and $\otimes$ denotes the Kronecker product,
see e.g. \cite[Section 4.2]{HJ2}. Minimality of the data implies
that the $2n$ points $x_1,~\ldots~,~x_n,~-x_1^*,~~\ldots~,~-x_n^*$
are distinct, so $\hat{X}$ is non-singular, see e.g.
\cite[Exercise 12]{HJ2}. This in turn implies that
$\mbox{\large\boldmath{X}}$ 
is nonsingular, see e.g.
\cite[Corollary 4.2.11]{HJ2}. Hence, the matrical coefficients
$C_o,~\ldots~,~C_{2n-1}$ in \eqref{P}
are unique and explicitly obtained.
\vskip 0.2cm

Finally, the structure guarantees that
\[
{C_k}^*=(-1)^kAC_kB\quad\quad\quad\quad\quad 
k=0,~\ldots~,~2n-1~.
\]
\qed
\vskip 0.2cm

In the sequel we focus on a special case of \eqref{InterpP}
where
\[
B=(A^*)^{-1}
\]
Namely $P(s)$ is so that,
\begin{equation}\label{InterpPstruct}
\begin{smallmatrix}
P^{\#}(s)&=&AP(s)(A^*)^{-1}&~&~&A\in\C^{m\times m}\quad
{\rm non-singular}\\~\\
P(x_j)&=&Y_j&~&~&j=1,~\ldots~,~n.
\end{smallmatrix}
\end{equation}

\section{constructing minimal degree symmetric neutral
polynomials}\label{sec:NeutPoly}
\setcounter{equation}{0}

We now address the problem of constructing, within a prescribed
family $\mathcal{F}$, polynomials $\Psi(s)$ of minimal degree,
vanishing at the given nodes $x_1,~\ldots~,~x_p\in\C$.
\vskip 0.2cm

Using \eqref{eq:RedData} let now $x_1,~\ldots~,~x_n$
be a resulting reduced set, see \eqref{eq:RedData}. Next, take
\begin{equation}\label{OriginalHatPsi}
\Psi(s):=\prod_{j=1}^{n}(x_j-s)M(x_j^*+s)
\quad\quad\quad\quad M\in\C^{m\times m},
\end{equation}
with $M$ parameter. Clearly, $\Psi(s)$ vanishes at the original
points $x_1,~\ldots~,~x_p$.
\vskip 0.2cm

Hence, taking $M$ to be: Hermitian or Skew-Hermitian
yields a minimal degree $\Psi$ in: ${\mathcal Even}$ or
${\mathcal Odd}$, respectively.
\vskip 0.2cm

Note that for $s=i\omega$ with $\omega$ real, $\Psi(s)$ in
\eqref{OriginalHatPsi} satisfies 
\begin{equation}\label{PsiOmega}
\Psi(i\omega)=M\prod_{j=1}^{n}|x_j-i\omega|^2.
\end{equation}

To guarantee
that the sets of polynomials described in \eqref{InterpPstruct}
and \eqref{OriginalHatPsi} indeed intersect, we need to further
restrict the parameter $M$ in \eqref{OriginalHatPsi} so that
the product $AM$, with $A$ from \eqref{InterpPstruct}, is
Hermitian, i.e.
\begin{equation}\label{eq:CommuteAM}
AM=(AM)^*.
\end{equation}

\section{interesting special cases}
\label{sec:SpecialCases}
\setcounter{equation}{0}

We here specialize the ``recipe" from Section
\ref{sec:introduction} to interesting classes
of polynomials.

\subsection{Even}\label{sec:EvenLagrange}

To guarantee feasibility
of the problem, one needs to substitute in \eqref{eq:feasible}
$A=B=I_m$ and thus obtain,
\[
\begin{smallmatrix}
x_j-x_k=0&~&\Longrightarrow&~&Y_j&=&Y_k\\~\\
x_j+x_k^*=0&~&\Longrightarrow&~&Y_j&=&Y_k^*
\end{smallmatrix}\quad\quad\quad p\geq k\geq j\geq 1.
\]
Taking in \eqref{InterpPstruct} $A=I_m$ results in $P(s)$ is in
${\mathcal Even}$, i.e. $P=P^{\#}$. 

\subsection{
$J-{\mathcal Even}$}\label{sec:even}
We shall find it convenient to denote by $J$ an arbitrary
$m\times m$ Hermitian involution, i.e.
\begin{equation}\label{J}
J=J^*=J^{-1}.
\end{equation}
Recall that $J$ is unitarily similar to
${\rm diag}\{-I_{\nu}~,~I_{m-\nu}\}$ with
$\nu\in[0,~m]$, see e.g. \cite[Theorem 4.1.5]{HJ1}.
\vskip 0.2cm

Taking in \eqref{InterpPstruct} $A=J$ with $J$ as in
\eqref{J}, results in $P(s)$ in $J-{\mathcal Even}$.
\vskip 0.2cm

In order to have the problem feasible one needs to assume that
the original data satisfies
\[
\begin{smallmatrix}
x_j+x_k^*=0&~&\Longrightarrow&~&Y_j&=&JY_k^*J
\end{smallmatrix}\quad\quad\quad p\geq k\geq j\geq 1.
\]
\vskip 0.2cm

For example if $m$ is even and $J=${\mbox{\tiny$\begin{pmatrix}
0&-I_{\frac{m}{2}}\\~\\ I_{\frac{m}{2}}&0\end{pmatrix}$}} then on
the imaginary axis $F(s)$ has Hamiltonian structure, see e.g.
\cite[Section I]{DMPP},
\cite[Problem 3.21, Eq. (6.3.3), Theorem 11.5.1]{GL} or
for a thorough treatment of this structure,
\cite[Section 7.2]{LR}.
\vskip 0.2cm

\subsection{Odd}\label{sec:Odd}

Taking in \eqref{InterpPstruct} $A=iI_m$ results in $P(s)$
in ${\mathcal Odd}$. In particular, in order to have
the problem feasible one needs to assume that the original data
satisfies
\[
\begin{smallmatrix}
x_j+x_k^*=0&~&\Longrightarrow&~&Y_j&=&-Y_k^*
\end{smallmatrix}\quad\quad\quad p\geq k\geq j\geq 1.
\]
In this case, the condition in \eqref{eq:CommuteAM} implies
that $M$ in \eqref{OriginalHatPsi} is skew-Hermitian.

\section{minimal degree interpolating Generalized Positive Even
polynomials}\label{sec:GPE}
\setcounter{equation}{0}

%
Let $x_1~,~\ldots~,~x_p\in\C$ and
$Y_1, \ldots~,~Y_p\in\C^{m\times m}$ be a feasible data set
i.e.
\[
\begin{smallmatrix}
x_j+x_k^*=0&~&\Longrightarrow&~&Y_j=Y_k&~&~&p\geq k>j\geq 1\\~\\
x_j\in{i}\R&~&\Longrightarrow&~&Y_j\in\overline{\mathbb P}_m
&~&~&j=1,~\ldots~,~p.
\end{smallmatrix}
\]
One searches~ {\em all} minimal degree interpolating
$\mathcal{GPE}$ polynomials, $F(s)$ in \eqref{eq:BasicInterp}.
To simplify presentation assume that the data is already
{\em reduced}, i.e.
\[
p=n.
\] 
Now take: (i) From Subsection \ref{sec:EvenLagrange} a minimal
degree interpolating $P\in{\mathcal Even}$. (ii) From
\eqref{OriginalHatPsi}
all minimal degree
$\Psi\in\mathcal{GPE}$ vanishing at the interpolation points:
\begin{equation}\label{HatPsiGPE}
\Psi(s):=\prod_{j=1}^{n}(x_j-s)M(x_j^*+s)
\quad\quad\quad\quad M\in\overline{\mathbb P}_m~.
\end{equation}
We now establish the minimality of the degree of 
$\Psi(s)$ in \eqref{HatPsiGPE}.

\begin{Pn}\label{Pn:GpeNeut}
Let $\Psi(s)$ be a $m\times m$-valued, minimal degree
$\mathcal{GPE}$ polynomial vanishing at a given
set of distinct points \mbox{$x_1,~\ldots~,~x_p\in\C$.}
Let also \mbox{$x_1,~\ldots~,~x_n\in
\C
$}
be a corresponding reduced set, see \eqref{eq:RedData}.
Then $\Psi(s)$ is of the form \eqref{HatPsiGPE}.
\end{Pn}

{\bf Proof :}\quad
From a given data points $x_1,~\ldots~,~x_p$ let us denote
\[
g(s):=\prod_{j=1}^{p}(x_j-s).
\]
Recall that every $\mathcal{GPE}$ function admits a factorization
of the form \eqref{SpectFact}. Thus, every $m\times m$-valued
$\mathcal{GPE}$ polynomial $\Psi(s)$ vanishing at
\mbox{$x_1,~\ldots~,~x_p$} is of the form
\[
\Psi(s)=g(s)\tilde{\Psi}(s)g^{\#}(s)
\quad\quad\quad\quad\tilde{\Psi}\in\mathcal{GPE}.
\]
Let now $x_1,~\ldots~,~x_n$ be a {\em reduced} subset of
the data points and then
\[
g_o(s):=\prod_{j=1}^{n}(x_j-s).
\]
Hence, without loss of generality one can take
\[
\Psi(s)=g_o(s)\tilde{\Psi}(s)g_o^{\#}(s)
\quad\quad\quad\quad\tilde{\Psi}\in\mathcal{GPE}.
\]
To guarantee minimality of the degree of $\tilde{\Psi}(s)$,
take 
it to be of degree zero, i.e.
\[
\tilde{\Psi}(s)\equiv M\in\overline{\mathbb P}_m~,
\]
so the claim is established.
\qed
\vskip 0.2cm

Without loss of generality we shall find it convenient to
normalize $M$ in \eqref{HatPsiGPE} so that
\[
\| M\|_2=1.
\]
Following the ``recipe" from Section \ref{sec:introduction},
the sought interpolating polynomial $F(s)$ is of the form
\begin{equation}\label{eq:Scheme}
F(s)=P(s)+\beta\Psi(s),
\end{equation}
with $P(s)$, $\Psi(s)$ from Propositions \ref{Pr:Lagrange}
\ref{Pn:GpeNeut}, respectively and $\beta$ a parameter.
\vskip 0.2cm

Assume that $P(s)$ in \eqref{P} is not in $\mathcal{GPE}$.
First, we fix in \eqref{HatPsiGPE} an arbitrary $M$ in
${\mathbb P}_m$ (i.e. non-singular). By construction, for all
$\beta\in\R$ in \eqref{eq:Scheme} is an interpolating
polynomial \eqref{eq:BasicInterp} in ${\mathcal Even}$.
\vskip 0.2cm

Now, on the one hand in \eqref{P} $P\in{\mathcal Even}$ is of
degree of (at most) $(2n-1)m$. On the other hand, in
\eqref{HatPsiGPE} $\Psi(s)$ is of degree $2nm$ and in
$\mathcal{GPE}$, a convex subset of ${\mathcal Even}$. 
Thus, there exists $\hat{\beta}>0$ so that
in \eqref{eq:InterpPolynom}, \eqref{eq:Scheme}
\[
\Psi\in\mathcal{GPE}\quad\quad\quad\forall\beta\geq\hat{\beta}.
\]
We next show that here $\hat{\beta}$ can be explicitly obtained.
\vskip 0.2cm

\begin{Pn}\label{Pn:HatBeta}
Assume that $P(s)$ in \eqref{P} is not in $\mathcal{GPE}$ and
that in \eqref{HatPsiGPE} $M\in{\mathbb P}_m$ is given. Then
$\hat{\beta}$ in \eqref{eq:InterpPolynom}, \eqref{eq:Scheme}
is given by
\begin{equation}\label{hatBeta}
\hat{\beta}=-\min\limits_{\begin{smallmatrix}\omega\in\R\\
\omega\not=-ix_j\end{smallmatrix}}
\left(\prod\limits_{j=1}^{n}|x_j-i\omega|^{-2}
\min\limits_{i=1~,~\ldots~,~m}~\lambda_i\left(
M^{-1}\sum\limits_{k=0}^{2n-1}i^kC_k\omega^k\right)\right)
\end{equation}
\end{Pn}

{\bf Proof }\quad
Indeed, note that $F(s)$ in \eqref{eq:InterpPolynom},
\eqref{eq:Scheme} is in
$\mathcal{GPE}$ if the following relations hold
\[
\begin{smallmatrix}
{\rm for~almost~all}~s\in{i}\R&F(s)&\in\overline{\mathbb P}_m
\\~\\
{\rm for~almost~all}~\omega\in\R&F(i\omega)&\in\overline{\mathbb P}_m
\\~\\
{\rm for~almost~all}~\omega\in\R&\left(\sum\limits_{k=0}^{2n-1}
C_k(i\omega)^k
+\beta\prod\limits_{j=1}^{n}|x_j-i\omega|^2M\right)
&\in\overline{\mathbb P}_m
\\~\\
{\rm for~almost~all}~\omega\in\R&\left(\sum\limits_{k=0}^{
2n-1}M^{-\frac{1}{2}}C_kM^{-\frac{1}{2}}{i\omega}^k
+\beta\prod\limits_{j=1}^{n}|x_j-i\omega|^2
\right) &\in\overline{\mathbb P}_m
\\~\\
{\rm for~almost~all}~\omega\in\R&\left(\frac{
\sum\limits_{k=0}^{2n-1}
M^{-\frac{1}{2}}C_kM^{-\frac{1}{2}}{i\omega}^k}
{\prod\limits_{j=1}^{n}|x_j-i\omega|^2}+\beta{I}_m\right)&
\in\overline{\mathbb P}_m
\end{smallmatrix}
\]
This can be written as,
\[
\begin{smallmatrix}
{\rm for~almost~all}~\omega\in\R\quad\min\limits_{i=1,~\ldots~,~m}
\lambda_i&\left(\frac{\sum\limits_{k=0}^{2n-1}
M^{-\frac{1}{2}}C_kM^{-\frac{1}{2}}{i\omega}^k}
{\prod\limits_{j=1}^{n}|x_j-i\omega|^2}+\beta{I}_m\right)&\geq 0
\\~\\
{\rm for~almost~all}~\omega\in\R\quad\min\limits_{i=1,~\ldots~,~m}
\lambda_i&\left(\frac{\sum\limits_{k=0}^{2n-1}
M^{-\frac{1}{2}}C_kM^{-\frac{1}{2}}{i\omega}^k}
{\prod\limits_{j=1}^{n}|x_j-i\omega|^2}\right)&\geq-\beta
\\~\\
{\rm for~almost~all}~\omega\in\R\quad\min\limits_{i=1,~\ldots~,~m}
\lambda_i&\left(\frac{M^{-1}\sum\limits_{k=0}^{2n-1}
C_k{i\omega}^k}
{\prod\limits_{j=1}^{n}|x_j-i\omega|^2}\right)&\geq-\beta
\\~\\
{\rm for~almost~all}~\omega\in\R
\left(\prod\limits_{j=1}^{n}|x_j-i\omega|^{-2}
\min\limits_{i=1,~\ldots~,~m}\lambda_i\right.&\left.\left(
M^{-1}\sum\limits_{k=0}^{2n-1}C_k{i\omega}^k
\right)\right)&\geq-\beta .
\end{smallmatrix}
\]
Thus, the claim is established.
\qed
\vskip 0.2cm

The search for minimum over all $\omega\in\R$ may in practice
be confined to a small interval as
\[
\lim\limits_{\omega\rightarrow~\pm\infty}
\left(
\frac{M^{-1}\sum\limits_{k=0}^{2n-1}i^kC_k\omega^k}
{\prod\limits_{j=1}^{n}
|x_j-i\omega|^2}\right)=0.
\]
\vskip 0.2cm

The above construction is illustrated by the following
example.
\vskip 0.2cm

\begin{Ex}\label{Ex:Motivation4}
{\rm

Consider the problem from
Examples \ref{Ex:Motivation1},
\ref{Ex:Motivation2} and \ref{Ex:Motivation3} of finding $F(s)$
so that
\[
\begin{smallmatrix}
F(s)&~&
\left(\begin{smallmatrix}1\\ 2\\ 3\end{smallmatrix}\right)&
\longrightarrow&
\left(\begin{smallmatrix}18\\ 75\\ 50\end{smallmatrix}\right).
\end{smallmatrix}
\]
From \eqref{P} we have that here
\[
\left(\begin{smallmatrix}
1&~~1&~1&1&1&1\\
1&~~2&~4&8&16&32\\
1&~~3&~9&27&81&243\\
1&-1&~1&-1&1&-1\\
1&-2&~4&-8&16&-32\\
1&-3&~9&-27&81&-243
\end{smallmatrix}\right)
\left(\begin{smallmatrix}C_o\\C_1\\C_2\\C_3\\C_4\\C_5
\end{smallmatrix}\right)=\left(\begin{smallmatrix}
18\\75\\50\\18\\75\\50\end{smallmatrix}\right)
\]
and thus
\[
P_4(s)=-3s^4+34s^2-13
\]
and from \eqref{HatPsiGPE} $\Psi(s)$ is as in \eqref{eq:ExHatPsi}.
Hence, for all $\beta\in\R$,
\begin{equation}\label{eq:TruePsib}
F_4(s)=P_4(s)+\beta\Psi(s)=-3s^4+34s^2-13+
\beta(1-s^2)(4-s^2)(9-s^2),
\end{equation}
is an interpolating polynomial in ${\mathcal Even}$. Next note that
$P_4(s)$ is not in $\mathcal{GPE}$,
\[
P_4(s)_{|_{s=i\omega}}=-(3\omega^4+34\omega^2+13)
\]
namely in fact $-P\in\mathcal{GPE}$.
Next, from \eqref{eq:TruePsib} one has that
\[
F_4(s)_{|_{s=i\omega}}=\left(P_4(s)+\beta\Psi(s)
\right)_{|_{s=i\omega}}=\left(\beta-\frac{1}{2}\right)(
\omega^6+14\omega^4+49\omega^2+36)
+\frac{1}{2}(\omega^2-1)^2(\omega^2+10).
\]
Thus, $F_4\in\mathcal{GPE}$ for all $\beta\geq\frac{1}{2}$.
Furthermore, on the boundary, i.e. for $\beta=\frac{1}{2}$
\[
\min\limits_{\omega\in\R}{\rm Re}~\left(F_4(s)_{|_{s=i\omega}}
\right)={F_4(s)}_{|_{s=\pm{i}}}=0.
\]
Thus indeed one obtains ~{\em all} interpolating $\mathcal{GPE}$
polynomials.
\vskip 0.2cm

Obviously, formally
applying \eqref{hatBeta} leads to the same conclusion.
\vskip 0.2cm

Recall that in Example \ref{Ex:Motivation3}, through
the simpler Lagrange approach we identified in \eqref{eq:ExPsib}
only a subset of interpolating functions with
$\beta\geq\frac{5}{4}~$.
}
\qed
\end{Ex}

\section{a refinement}\label{sec:refine}
\setcounter{equation}{0}

We here address a case enabling us to refine the construction
in the previous section by allowing in \eqref{HatPsiGPE}
$M\in\overline{\mathbb P}_m$ (singular) and consequently obtain
lower degree interpolating polynomials $F(s)$.
\vskip 0.2cm

Recall that
in subsection \ref{sec:EvenLagrange} 
$P(s)$ obtained, is an interpolating polynomial in
${\mathcal Even}$. Thus on $i\R$ it is Hermitian. Assume that
there exists a (constant) nonsingular
matrix $T$ so that
\begin{equation}\label{DiagP}
TP(s)T^*={\rm diag}\{P_r(s),~P_{m-r}(s)\}\quad\quad\quad\quad
P_{m-r}\in\mathcal{GPE},~r\in[0, m-1].
\end{equation}
(The case addressed in the previous section corresponds to
$r=m$.)
\vskip 0.2cm

If $r=0$, $P(s)$ in \eqref{P} is already in $\mathcal{GPE}$, it
is a minimal degree interpolating polynomial and one can now
employ \eqref{eq:Scheme} with arbitrary $\Psi(s)$
from \eqref{HatPsiGPE} and arbitrary $\beta\geq 0$. This is
illustrated next.
\vskip 0.2cm

\begin{Ex}\label{Ex:Psia}

{\rm
Consider the scalar problem 
of finding ~{\em all}~ $\mathcal{GPE}$ polynomials $F(s)$ so that
\[
\begin{smallmatrix}
F(s)&~&\left(\begin{smallmatrix}1\\ 2\\ 3\end{smallmatrix}\right)&
\longrightarrow&
\left(\begin{smallmatrix}~~~4\\~~~1\\-4\end{smallmatrix}\right).
\end{smallmatrix}
\]
From \eqref{P} we have that here
\[
\left(\begin{smallmatrix}
1&~~1&~1&1&1&1\\
1&~~2&~4&8&16&32\\
1&~~3&~9&27&81&243\\
1&-1&~1&-1&1&-1\\
1&-2&~4&-8&16&-32\\
1&-3&~9&-27&81&-243
\end{smallmatrix}\right)
\left(\begin{smallmatrix}C_o\\C_1\\C_2\\C_3\\C_4\\C_5
\end{smallmatrix}\right)=\left(\begin{smallmatrix}~~
4\\~~1\\-4\\~~4\\~~1\\-4\end{smallmatrix}\right)
\]
and thus,
\[
P(s)=-s^2+5
\]
and from \eqref{HatPsiGPE}
\begin{equation}\label{eq:ExHatPsi}
\Psi(s)=(1-s^2)(4-s^2)(9-s^2).
\end{equation}
As we here have $P(i\omega)=\omega^2+5$ for $\omega\in\R$,
$P\in\mathcal{GPE}$ and all interpolating polynomials are given by
\begin{equation}\label{eq:TruePsia}
F(s)=P(s)+\beta\Psi(s)=
-s^2+5+\beta(1-s^2)(4-s^2)(9-s^2)\quad\quad\quad\beta\geq 0.
\end{equation}
It is of interest to mention that employing the simpler Lagrange
approach described in Example \ref{Ex:Motivation3} would identify
the subset of interpolating polynomials $F(s)$ in \eqref{eq:TruePsia}
with $\beta\geq\frac{1}{36}$. In particular, it would have failed to
find the minimal degree interpolating
$\mathcal{GPE}$ polynomial corresponding to $\beta=0$.
}\qed
\end{Ex}
\vskip 0.2cm

Consider now the case where in \eqref{DiagP}
\[
r\in[1, m-1].
\]
Then, conforming to $P(s)$ in \eqref{DiagP}, $\Psi(s)$ in
\eqref{HatPsiGPE} can be constructed, with the same $T$, so
that
\[
T\Psi(s)T^*=\prod_{j=1}^{n}(x_j-s){\rm diag}\{M_r,~M_{m-r}\}
(x_j^*+s)\quad\quad\quad\quad M_r\in{\mathbb P}_r~,
\]
where $M_{m-r}\in\overline{\mathbb P}_{m-r}$ is arbitrary,
including zero.
\vskip 0.2cm

One can now proceed as before. Namely, fix $M_r\in{\mathbb P}_r$,
where \mbox{$\| M_r\|_2=1$}. Then
$\hat{\beta}$ in 
is given by
\[
\hat{\beta}=-\min\limits_{\begin{smallmatrix}\omega\in\R\\
\omega\not=ix_j\end{smallmatrix}}\left(
\prod\limits_{j=1}^{n}|x_j-i\omega|^{-2}
\min\limits_{i=1~,~\ldots~,~r}~\lambda_i\left(M_r^{-1}P_r(i\omega)
\right)\right).
\]

\begin{Ex}\label{Ex:Psic}
{\rm
Consider the following two dimensional problem: Find all $F(s)$ in
$\mathcal{GPE}$ so that
\[
\begin{smallmatrix}
F(1)&=&{\rm diag}\{-35,&~~9\}\\~\\
F(2)&=&{\rm diag}\{-20,&~~0\}\\~\\
F(3)&=&{\rm diag}\{45,&~~25\}
\end{smallmatrix}
\]
From \eqref{P} we have that here
\[
\left(\begin{smallmatrix}
I_2&~~I_2&~I_2&I_2&I_2&I_2\\
I_2&~~2I_2&~4I_2&8I_2&16I_2&32I_2\\
I_2&~~3I_2&~9I_2&27I_2&81I_2&243I_2\\
I_2&-I_2&~I_2&-I_2&I_2&-I_2\\
I_2&-2I_2&~4I_2&-8I_2&16I_2&-32I_2\\
I_2&-3I_2&~9I_2&-27I_2&81I_2&-243I_2
\end{smallmatrix}\right)
\left(\begin{smallmatrix}C_o\\C_1\\C_2\\C_3\\C_4\\C_5
\end{smallmatrix}\right)=\left(\begin{smallmatrix}
{\rm diag}\{-35&~~9\}\\
{\rm diag}\{-20&~~0\}\\
{\rm diag}\{45&~~25\}\\
{\rm diag}\{-35&~~9\}\\
{\rm diag}\{-20&~~0\}\\
{\rm diag}\{45&~~25\}\end{smallmatrix}\right).
\]
Thus $C_1$, $C_3$, $C_5$ vanish and
\[
C_o={\rm diag}\{-36\quad 16\}\quad\quad\quad
C_2={\rm diag}\{0\quad -8\}\quad\quad\quad
C_4=I_2~.
\]
Namely,
\[
\begin{smallmatrix}
P(s)&=&I_2s^4+{\rm diag}\{0,~~-8\}s^2+
{\rm diag}\{-36,~~16\}\\~\\~
&=&{\rm diag}\{s^4-36,~~(4-s^2)^2\},\end{smallmatrix}
\]
is a minimal degree interpolating polynomial in ${\mathcal Even}$.
Using \eqref{eq:Deg} the McMillan degree is 10.
\vskip 0.2cm

Note now that $P(s)$ is of the form of \eqref{DiagP} with
$T=I_2$ and $r=1$.
\vskip 0.2cm

Now from \eqref{HatPsiGPE}  one can construct $\Psi(s)$ with
$M\in\overline{\mathbb P}_2$. Indeed taking
\[
\Psi(s)=(1-s^2)(4-s^2)(9-s^2)\left(
{\rm diag}\{\beta,~~0\}+\Delta\right)
\quad\quad\quad\quad
\quad\Delta\in\overline{\mathbb P}_2~,
\]
guarantees that
\[
F(s)=P(s)+\Psi(s)
\]
is a $\mathcal{GPE}$ interpolating polynomial for all
\[
\beta\geq\hat{\beta}=1.
\]
Moreover for $\Delta=0$ one obtains
\[
\begin{smallmatrix}
F(s)&=&{\rm diag}\{s^4-36+
\beta(1-s^2)(4-s^2)(9-s^2),~~(4-s^2)^2\}\\~\\~
&=&
{\rm diag}\{
-\beta{s}^6+(14\beta+1)s^4-49\beta{s}^2+36(\beta-1),~~(4-s^2)^2\}
\\~\\~&=&
s^6{\rm diag}\{-\beta,~~0\}
+s^4{\rm diag}\{14\beta+1,~~1\}
+s^2{\rm diag}\{-49\beta,~~-8\}
+{\rm diag}\{36(\beta-1),~~16\}
\end{smallmatrix}
\]
where the McMillan degree is only 12 (see\eqref{eq:Deg}). For
$\beta\geq 1$ indeed $F\in\mathcal{GPE}$.
}\qed
\end{Ex}

\section{a non-convex set: $\nu$-generalized positive even
polynomials}\label{sec:J-GPE}
\setcounter{equation}{0}

We here address the the interpolation problem \eqref{eq:BasicInterp}
where $\mathcal{F}$ is not convex. We modify the above ``recipe"
accordingly. As a test case, we take the set $\mathcal{{\nu}GPE}$
described in \eqref{JSpectFact}. First, we show that this set
is indeed not convex.
\vskip 0.2cm

As already mentioned, $\pm\mathcal{GPE}$ are sub{\bf cic}s of
${\mathcal Even}$, see e.g.
\cite[Section 5]{AL3}. However, for $\nu\in[1,~m-1]$ the
set $\mathcal{{\nu}GPE}$ in \eqref{JSpectFact} is
an invertible cone, but not convex. This can be illustrated
even by constant $2\times 2$ matrices. Indeed if one takes
\[
\begin{smallmatrix}
{\rm diag}\{-1,~4\}&=&A&=&R_a{\rm diag}\{-1,~1\}R_a^*&~&R_a&
=&{\rm diag}\{1,~2\}\\~\\
{\rm diag}\{4,~-1\}&=&B&=&R_b{\rm diag}\{-1,~1\}R_b^*&~&
R_b&=&\left(\begin{smallmatrix}0&2\\1&0\end{smallmatrix}\right)
\end{smallmatrix}
\]
then
$A+B=3I_2$ is not a $\mathcal{{\nu}GPE}$ function
of the form of \eqref{JSpectFact}.
\vskip 0.2cm

The fact that the set $\mathcal{{\nu}GPE}$ is not convex
does not allow us to employ the Lagrange approach from
Section \ref{sec:InterpolationBackground}. However, we can
adapt the recipe from Section \ref{sec:introduction}.
Here are the details.
\vskip 0.2cm

Given a reduced data set \mbox{$x_1~,~\ldots~,~x_n\in\C$,}
the corresponding \mbox{$Y_1~,~\ldots~,~Y_n\in\C^{m\times m}$}
and $\nu$, $\nu\in[1,~m-1]$.  Find (all) $F(s)$, low degree
${m\times m}$-valued interpolating polynomials
$F~:~\C~\rightarrow~\C^{m\times m}$, i.e.
\[
\begin{smallmatrix}
F(s)&=&
G(s){\rm diag}\{-I_{\nu}~,~I_{m-\nu}\}G^{\#}(s)
&~&~&\nu\in[1,~m-1]\quad G(s)~{\rm polynomial}\\~\\
F(x_j)&=&Y_j&~&~&j=1,~\ldots~,~n.
\end{smallmatrix}
\]
In order to have the problem feasible one needs to assume that
the data satisfies
\[
\begin{smallmatrix}
x_j+x_k^*=0&~&\Longrightarrow&~&Y_j=Y_k&~&~&n\geq k>j\geq 1\\~\\
x_j\in{i}\R&~&\Longrightarrow&~&Y_j=T_j{\rm
diag}\{-I_{\nu}~,~I_{m-\nu}\}T_j^*&~&~&j=1,~\ldots~,~n,
\end{smallmatrix}
\]
for some\begin{footnote}{Recall, if $\nu=0$ or $\nu=l$ we are
essentially back to the $\mathcal{GPE}$ case of subsection
\ref{sec:GPE}.}\end{footnote} $\nu\in[1,~m-1]$ and some
$T_j\in\C^{m\times m}$.
\vskip 0.2cm

Substituting in subsection \ref{sec:even} $A=J=I_m$ one obtains
from \eqref{P} $P(s)$, the minimal degree interpolating
polynomial in ${\mathcal Even}$.
\vskip 0.2cm

To construct the neutral polynomials $\Psi(s)$ substitute in
\eqref{OriginalHatPsi}
\mbox{$M=R{\rm diag}\{-I_{\nu}~,~I_{m-\nu}\}R^*$} with
$\nu\in[1,~m-1]$ and $x_j$ as above, to obtain
\begin{equation}\label{HatPsiJGPE}
\Psi(s):=\prod_{j=1}^{n}(x_j-s)\left(R{\rm diag}\{-I_{\nu}~,~
I_{m-\nu}\}R^*\right)(x_j^*+s)
\quad\quad\quad\quad R\in\C^{m\times m},
\end{equation}
with $R$ parameter.
\vskip 0.2cm

By construction, with the above $P(s)$ and $\Psi(s)$, for all
$\beta\in\R$
\begin{equation}\label{eq:PsiJGPE}
F(s)=P(s)+\beta\Psi(s),
\end{equation}
is an interpolating polynomial \eqref{eq:BasicInterp}
in ${\mathcal Even}$. However, as we already remarked,
the set $\mathcal{{\nu}GPE}$ not convex. Thus one needs to justify
the existence of $\hat{\beta}$ so that $F(s)$ in
\eqref{eq:PsiJGPE} is a $\mathcal{{\nu}GPE}$ interpolating
polynomial for all $\beta>\hat{\beta}$.
\vskip 0.2cm

The idea relies on the following fact, here formulated in the
framework of matrix theory.

\begin{La}
Let $A_p, A_{\Psi}\in\C^{m\times m}$ Hermitian matrices.
Assume that $A_{\Psi}$ has $\pi$ and $m-\pi$ eigenvalues in
$\C_+$ and $\C_-$, respectively ($A_{\Psi}$ is nonsingular).
Then, for all
\[
\beta>\|A_p\| \|A_{\Psi}^{-1}\|
\]
the matrix $(A_p+\beta{A}_{\Psi})$ has $\pi$ eigenvalues in
$\C_+$ and $m-\pi$ eigenvalues in $\C_-$.
\end{La}

This may be deduced in several ways\begin{footnote}{In operator
theory this is formulated as having $A_p+{A}_{\Psi}$ invertible
whenever $A_{\Psi}$ is invertible and
\mbox{$\|A_{\Psi}^{-1}\|^{-1}>\|A_p\|$}, see e.g.
\cite[Theorem 10.20]{Ru}.}\end{footnote}, e.g. Weyl's Theorem
\cite[Theorem 4.3.1]{HJ1}. A detailed proof in our
framework, is given in Proposition \ref{Pn:EstimateHatBeta}
below where we estimate $\hat{\beta}$.

\begin{Pn}\label{Pn:EstimateHatBeta}
Let $F(s)$ be as in \eqref{eq:PsiJGPE} with $P(s)$ from \eqref{P}
and $\Psi(s)$ as in \eqref{HatPsiJGPE} where $R$ is given and
nonsingular.

Then, in \eqref{eq:PsiJGPE} $F\in\mathcal{{\nu}GPE}$
for all $\beta>\hat{\beta}$ where
\begin{equation}\label{eq:HatBetaJGPE}
\hat{\beta}=\max\limits_{\omega\in\R}
\left(
\prod\limits_{j=1}^{n}|x_j-i\omega|^{-2}\left\|
\sum\limits_{k=0}^{2n-1}R^{-1}C_k(R^*)^{-1}
(i\omega)^k\right\|\right).
\end{equation}
If $x_j=i\omega_j$ with $\omega_j\in\R$ is an interpolation node,
a whole neighborhood of $\omega_j$ is excluded from the above
$\max\limits_{\omega\in\R}$.
\end{Pn}

{\bf Proof}\quad
Indeed, note that $F(s)$ in \eqref{eq:PsiJGPE} is in
$\mathcal{{\nu}GPE}$ if there exist a polynomial $G(s)$
so that
\[
\begin{smallmatrix}
{\rm for~almost~all}~s\in{i}\R&\left(P(s)+\beta\Psi(s)\right)
&=G(s){\rm diag}\{-I_{\nu}~,~I_{m-\nu}\}G^{\#}(s)
\\~\\
{\rm for~almost~all}~s\in{i}\R&\left(\sum\limits_{k=0}^{2n-1}
C_ks^{k}+\beta\prod\limits\limits_{j=1}^{n}(x_j-s)R
{\rm diag}\{-I_{\nu}~,~I_{m-\nu}\}R^*(x_j^*+s)\right)
&=G(s){\rm diag}\{-I_{\nu}~,~I_{m-\nu}\}G^{\#}(s).
\end{smallmatrix}
\]
If $s=x_j\in{i}\R$ is an interpolation node, then
\[
\sum\limits_{k=0}^{2n-1}{C_ks^{k}}_{|_{s=x_j}}=Y_j=T_j{\rm
diag}\{-I_{\nu}~,~I_{m-\nu}\}T_j^*
\]
and
\[
\prod\limits\limits_{j=1}^{n}(x_j-s)R{\rm
diag}\{-I_{\nu}~,~I_{m-\nu}\}R^*(x_j^*+s)=0,
\]
so the condition is satisfied in a neighborhood of $x_j$. Hence,
assume hereafter that $s\in{i}\R$ is out of a neighborhood of
interpolation points. Thus, the above condition may be written
as having a rational function $\tilde{G}(s)$ so that
\[
\begin{smallmatrix}
{\rm for~almost~all}~s\in{i}\R&\left(\frac
{\sum\limits_{k=0}^{2n-1}R^{-1}C_k(R^*)^{-1}
s^k}{\prod\limits_{j=1}^{n}(x_j-s)(x_j^*+s)}
+\beta{\rm diag}\{-I_{\nu}~,~I_{m-\nu}\}\right)
&=\tilde{G}(s){\rm diag}\{-I_{\nu}~,~I_{m-\nu}\}\tilde{G}^{\#}(s)
\end{smallmatrix}
\]
Now, this in turn is implied by,
\[
\begin{smallmatrix}
{\rm for~almost~all}~s\in{i}\R&~&~&\beta\geq
\left\|\frac{\sum\limits_{k=0}^{2n-1}R^{-1}C_k(R^*)^{-1}
s^k}{\prod\limits_{j=1}^{n}(x_j-s)(x_j^*+s)}\right\|.
\end{smallmatrix}
\]
Substituting $s=i\omega$, $\omega\in\R$,
establishes the claim.
\qed
\vskip 0.2cm

Note that the value of $\hat{\beta}$ in \eqref{eq:HatBetaJGPE}
depends on the choices of
$R$ in \eqref{HatPsiJGPE} and of the norm in
\eqref{eq:HatBetaJGPE}.

\section{future research}
\label{sec:FutureResearch}
\setcounter{equation}{0}

As it is often the case, this study opens the door for future
research problem. We here mention a sample of them.

\begin{itemize}
\item[(i)~~~]{}Recall that in \eqref{SpectFact} we mentioned
that $F\in\mathcal{GPE}$ if and only if it admits a
factorization of the form $F(s)=G(s)G^{\#}(s)$.

\noindent
In Sections \ref{sec:Lagrange} through \ref{sec:refine}
we presented a ``factorization free" recipe for obtaining all
minimal degree interpolating $\mathcal{GPE}$ polynomials.

\noindent
It is now of interest, for a given minimal degree
interpolating $\mathcal{GPE}$ polynomial $F(s)$, to explore
properties of its (pseudo) spectral factors $G(s)$, see e.g.
\cite[Section 5.2]{AV}, \cite{BGKR} and
\cite[Section 19.3]{LR}. 
\vskip 0.2cm

\item[(ii)~~]{}An idea we have used throughout the work is as
follows. If $P(s)$ maps \mbox{$x_1,~\ldots~,~x_n$} to
\mbox{$Y_1,~\ldots~Y_n$} and $\Psi(s)$ vanishes at
\mbox{$x_1,~\ldots~,~x_n$} taking
\begin{equation}\label{eq:Basic}
F(s)=P(s)+\beta\Psi(s)
\end{equation}
yields $F(s)$ mapping \mbox{$x_1,~\ldots~,~x_n$} to
\mbox{$Y_1,~\ldots~Y_n$}, for all $\beta\in\C$.
\vskip 0.2cm

\noindent
Let now $\mathcal{F}_{\Psi}$, $\mathcal{F}_{P}$ be two families
of functions where $\mathcal{F}_{\Psi}$ is a (convex) subcone
of $\mathcal{F}_{P}$. Assuming $P\in\mathcal{F}_{P}$ and
$\Psi\in\mathcal{F}_{\Psi}$, for $\beta\geq 0$ ``sufficiently
large" $F(s)$ in \eqref{eq:Basic} is an interpolating function
within the family $\mathcal{F}_{\Psi}$.
\vskip 0.2cm

\noindent
In \cite{ABL0} we adapt this idea to interpolation by scalar
rational positive functions (as in the classical Nevanlinna-Pick
Interpolation) but where the nodes are within $\C_-$. It is shown
that there always exist interpolating functions of degree equal
to the number of nodes. Moreover, an easy-to-compute recipe of
constructing these functions is introduced.
\vskip 0.2cm

\item[(iii)~]{}Adapt the recipe in Section \ref{sec:introduction}
for interpolation by $\mathcal{GP}$ (not necessarily even)
polynomials. The framework of \eqref{eq:Basic} still holds
with $P(s)$ mapping $x_j$ to $Y_j$ and $\Psi(s)$ a minimal
degree $\mathcal{GP}$ polynomial vanishing at $x_j$.
However, there are two basic differences:

\noindent
(a) There is no restriction on the structure of $P(s)$.
Namely, it is no longer necessary to double its degree.

\noindent
(b) $\Psi(s)$ is no longer of the form of
\eqref{OriginalHatPsi}. For example, for two points
$x_1,x_2\in\C_+$
\[
\Psi(s)=(x_1-s)(x_2-s)(\theta{x}_1^*+(1-\theta)x_2^*+s)M
\]
where the scalar $\theta$, $\theta\in[0,~1]$ and the matrix
$M$, $(M+M^*)\in\overline{\mathbb P}$ are parameters.
\vskip 0.2cm

\item[(iv)~]{}Adapt the interpolation scheme of this work to
cope with ~{\em tangential}~ structured matrix-valued
polynomial interpolation. Namely substitute \eqref{eq:BasicInterp}
by \eqref{eq:RightTangential} and \eqref{eq:LeftTangential}.
\vskip 0.2cm

\item[(v)~~]{}One may be interested in interpolation by polynomials
whose symmetry on the imaginary axis is of a group type, e.g.
unitary, $J$-unitary, contraction or $J$-contraction, in the
spirit of e.g. \cite{AlGo1}, \cite{GR}. It is of interest to
solve the interpolation problem in \eqref{eq:BasicInterp}
where $\mathcal{F}$ is a family of matrix valued polynomials
with this symmetry.
\vskip 0.2cm

\item[(vi)~]{}In Section \ref{sec:introduction} we pointed out
that the set of $\mathcal{GPE}$ ~{\em functions}~ is a
sub{\bf cic} (Convex Invertible Cone) of ${\mathcal Even}$. If
one focuses, as in this work, on the respective subsets of
$\mathcal{GPE}$ and ${\mathcal Even}$~ {\em polynomials},~ 
invertibility is no longer relevant. Hence, out of the {\bf cic}
structure, it is only the Convex Cone part that can be used
(as we indeed did).

\noindent
In \eqref{eq:reverse} we recalled the notion of ~{\em reverse}~
of a polynomial. It is then easy to verify that if a polynomial
is in $\mathcal{GPE}$, so is the corresponding reverse polynomial,
\eqref{eq:reverse}.
In fact, it turns out that the set of $\mathcal{GPE}$ polynomials
is a sub-Convex Reversible Cone of ${\mathcal Even}$ polynomials.

\noindent
It is of interest to explore the Convex Reversible Cone structure
of $\mathcal{GPE}$ polynomials and then to try employ it to
interpolation.
\end{itemize}


\begin{thebibliography}{ZZ}

\bibitem{AADLW}D. Alpay, T.Ya. Azizov, A. Dijksma
H. Langer and G. Wanjala \newblock ``A Basic Interpolation
Problem for Generalized Schur Functions and Coisometric
Realizations", {\em {Operator {T}heory: {A}dvances and
{A}pplications}},Vol. 143, pp. 39-76,
\newblock Birkh{\" a}user Verlag, Basel, 2003.

\bibitem{ABD}D. Alpay, V. Bolotnikov and A. Dijksma,
``On the {N}evanlinna-{P}ick Interpolation Problem for
Generalized {S}tieltjes Functions", {\em Integ. Eq. \& Op.
Theory}, Vol. 30, pp. 379-408, 1998.

\bibitem{ABL0}D.~Alpay, V. Bolotnikov and I. Lewkowicz,~
``Wrong Side Interpolation by Positive Function",~ a preprint.

\bibitem{ABL1}D.~Alpay, V. Bolotnikov and Ph. Loubaton,~
``An Interpolation Problem with Symmetry and Related
Questions",~ {\em Zeitschrift f{\" u}r Analysis \&
Anwendungen},~
Vol. 15, pp. 19-29, 1996.

\bibitem{ABL2}D.~Alpay, V. Bolotnikov and Ph. Loubaton,~
``One Two-Sided Residue Interpolation for
Matrix-Valued $H_2$ Functions with Symmetries",~
{\em J. of Math. Anal. \& Appl.},~ Vol. 200, pp. 76-105,
1996.

\bibitem{ADLRS}D. Alpay, A. Dijksm, H. Langer, S. Reich and
D. Shoikhet, ``Boundary Interpolation and Rigidity for
Generalized Nevanlinna functions", {\em Math. Nach.},
Vol. 283, No. 3, pp. 335-364, 2010.

\bibitem{ADLS}D.~Alpay, A.~Dijksma, H. Langer and Y.
Shondin, \newblock ``The Schur Transformation for
Generalized Nevanlinna Functions: Interpolation and
Self-adjoint Operator Realizations", {\em Compl. Anal.
oper. theory}, Vol. 1, pp. 169-210, newblock
Birkh{\" a}user Verlag, Basel, 2007.

\bibitem{AlGo1} D.~Alpay and I.~Gohberg.
\newblock ``Unitary Rational Matrix Functions"
\newblock In I.~Gohberg, editor, {\em Topics in interpolation
theory of rational matrix-valued functions},
{\em Operator Theory: Advances and Applications},
Vol. 33, pp. 175--222. Birkh{\" a}user Verlag, Basel, 1988.

\bibitem{AL1}D.~Alpay and I. Lewkowicz, ~{\em An
easy-to-compute factorization of rational generalized
positive functions}, {\em Sys. Cont. Lett.}~ Vol. 59,
pp. 517-521, 2010.

\bibitem{AL2} D. Alpay and I. Lewkowicz,~ ``The
Positive Real Lemma and Construction of all Realizations
of Generalized Positive Rational Functions",~
{\em Sys. Cont Lett},~ Vol. 60, pp. 985-993, 2011.

\bibitem{AL3}D.~Alpay and I. Lewkowicz, ~``Convex cones of
generalized positive rational functions and the
Nevanlinna-Pick interpolation," to appear in
{\em Lin. Alg \& Appl.}.

Available at

\url{http://arxiv.org/abs/1010.0546}.

%
%

\bibitem{AmDe}A. Amirshadyan, and V.A. Derkach,
``Interpolation in Generalized {N}evanlinna and {S}tieltjes
Classes",~ {\em J. Operator Theory}, Vol. 42, No. 1, pp.
145-188, 1999.

\bibitem{ACL}A. Amiraslani, R.M. Coreless and P. Lancaster,
``Linearization of matrix polynomials expressed in
polynomial bases",~ {\em IMA J.  Num. Anal.},~ Vol. 29,
pp. 141-157, 2009.

\bibitem{AM}B.D.O. Anderson and J. B. Moore,
``Algebraic Structure of Generalized Positive
Real Matrices", {\em SIAM J. Control}, Vol. 6,
pp. 615-624, 1968.

\bibitem{AV}B.D.O. Anderson and S.  Vongpanitlerd, ~{\em
Networks Analysis and Synthesis, A Modern Systems Theory
Approach},~ Prentice-Hall, New Jersey, 1973.

\bibitem{ABKW}A.C. Antoulas, J.A. Ball, J. Kang and
J.C. Willems, ~``On the Solution of the Minimal Rational
Interpolation Problem",~ {\em Lin. Alg. \& Appl.}, ~
Vol. 137, pp. 511-573, 1990.

\bibitem{Ball}J. A. Ball, ``Interpolation Problems and
Loewner Types for Meromorphic Matrix Functions"
{\em Integ. Eq. \& Op. Theory}, Vol. 6, pp. 804-840, 1983.

\bibitem{BGR}J.A. Ball, I. Gohberg and L. Rodman,~
{\em Interpolation of Rational Matrix Functions},~
Vol.~44 of {\em {Operator {T}heory: {A}dvances and
{A}pplications}}, Birkh{\" a}user Verlag, Basel, 1990.

\bibitem{BH}J.A. Ball and J.W. Helton, ``Interpolation
problems of Pick-Nevanlinna and Loewner types for
meromorphic matrix functions: parameterization of the
set of all solutions", {\em Integ. Eq. \& Op. Theory},
Vol. 9, pp. 155-203, 1986.

\bibitem{BK}J.A. Ball and J. Kang,~ ``Matrix Polynomial
Solutions of Tangential Lagrange-Sylvester Interpolation
Conditions of Low McMillan Degree",~ {\em Lin. Alg.
\& Appl.}, ~Vol. 138, pp. 699-746, 1990.

\bibitem{BGKR}H. Bart, I. Gohberg, M.A. Kaashoek and
A.C.M. Ran,~ {\em A State Space Approach to Canonical
Factorization with Applications},~ {\em { Operator
{T}heory: {A}dvances and {A}pplications}}, Vol. 200,
Birkh{\" a}user Verlag, Basel, 2010.

\bibitem{Be}V. Belevich, {\em Classical Network Theory},
Holden Day, San-Francisco, 1968.

\bibitem{BT}J-P. Berrut and L.N. Trefethen,~ ``Barycentric
Lagrange Interpolation",~ {\em SIAM Review},~ Vol. 46,
pp. 501-517, 2004.

\bibitem{bol_oam}V. Bolotnikov, "Boundary rigidity for Some
Classes of meromorphic functions", {\em Oper.  Matrices},
Vol. 3, pp. 283--301, 2009.

\bibitem{bol_oam2}V. Bolotnikov, "A Multi-Point Degenerate
Interpolation Problem for Generalized Schur Functions",
{\em Oper. Matrices}, Vol. 4 (2010) pp. 151--191.

\bibitem{CL1}N. Cohen ~and ~I. Lewkowicz,~ ``Convex~
Invertible~ Cones ~and~ the~ Lyapunov~ Equation",~
{\em Lin. Alg. \& Appl.}, ~
Vol. 250, pp. 105-131, 1997.

\bibitem{CL2}N. Cohen ~and ~I. Lewkowicz,~ ``Convex~
Invertible~ Cones ~of~ State ~Space~ Systems",~
{\it Mathematics of Control
Signals and Systems},~ Vol. 10, pp. 265-285, 1997.

\bibitem{CL4}N. Cohen and I. Lewkowicz, ~``Convex~
Invertible~ Cones and Positive Real Analytic Functions",~
{\it Lin. Alg. \& Appl.}, Vol. 425, pp. 797-813, 2007.

\bibitem{CL5}N. Cohen and I. Lewkowicz,~ ``The Lyapunov
order for real matrices", {\it Lin. Alg. \& Appl.}, Vol.
430, pp. 1489-1866, 2009.

\bibitem{DGK}Ph. Delsarte, Y. Genin and Y. Kamp.
``Pseudo-Carath\'{e}odory functions and Hermitian
Toeplitz Matrices'', {\em Philips J. of Research},
Vol. 41, pp. 1-54, 1986.

\bibitem{CZ}Recommended Publications on Polynomial Methods
in Control. Available at:

\url{http://www.polyx.cz/publications.html}

\bibitem{DZ}M. S. Derevyagin and A.S. Zhedanov,~ ``An
Operator Approach to Multipoint Pad\'{e} Approximations",~
{\em J. Approx. The}, Vol. 157, pp. 70-88, 2009.

\bibitem{DHdS}V.A. Derkach, S. Hassi and H. de-Snoo,
``Operator models associated with Kac subclasses of generalized
Nevanlinna functions'', {\em Meth. Funct. Anal. \&
Topology},~ Vol. 5, pp. 65-87, 1999.

\bibitem{DD}V.A. Derkach and H. Dym, ``Bitangential
Interpolation in Generalized Schur Class'',~
{\em Complex Anal. \& Op. Theory},~ Vol. 4, pp. 761-765,
2010.

\bibitem{DMPP}L. Dieci, B. Morini, A. Papini and A. Pasquali
``On Real Logarithm of Nearby Matrices and Structured
Matrix Interpolation", ~{\em Applied Numerical Mathematics},~
Vol. 29, pp. 145-165, 1999.

\bibitem{DLLS1}A. Dijksma, H. Langer, A. Luger
and Yu. Shondin, ``A factorization result for generalized
Nevanlinna functions of class $~{\mathcal N}_{\kappa}$'', {\em
Integ. Eq. \& Op. Theory}, Vol. 36, pp. 121-124, 2000.

%
\bibitem{Fu2}P.A. Fuhrmann,~ {\em A Polynomial Approach to
Linear ~Algebra},~ Springer 1996.

\bibitem{Fu3}P.A. Fuhrmann, ``On Tangential Matrix
Interpolation", {\em  Lin. Alg. \& Appl.}, Vol. 433,
pp. 2018-2059, 2010.

\bibitem{GR}I. Gohberg and I. Rubinstein, ``Proper
Contractions and their Unitary Minimal Completions",~
{\em { Operator {T}heory: {A}dvances and {A}pplications}},
Vol. 33, pp. 223-247, Birkh{\" a}user Verlag, Basel, 1988.

\bibitem{GL}M. Green and D. Limebeer, ~{\em Linear Robust
Control},~ Prentice-Hall, 1995.

%
\bibitem{HS}M. Henrion and M. \u{S}ebek,~ ``Symmetric
Matrix Polynomial Equation: Interpolation Results"
{\em Automatica},~ Vol. 34, pp. 811-824, 1998.

\bibitem{Hi}N.J. Higham,~ {\em Functions of Matrices- Theory and
Computation},~ SIAM, 2008.

\bibitem{HJ1}R.A. Horn and C.R. Johnson, ~{\em Matrix
Analysis},~ Cambridge University Press, 1985.

\bibitem{HJ2}R.A. Horn and C.R. Johnson, ~{\em Topics in
Matrix Analysis},~ Cambridge University Press, 1991.

%
%
\bibitem{La}J.L. Lagrange, ``Le\c{c}ons \'{E}l\'{e}mentaires
sur les Math\'{e}matiques Donn\'{e}es \`{a} l'\'{E}cole Normale",
in Oeuvre de Lagrange, J-A Serret Ed., Paris France,
Gauthier-Villars, Vol. 7, pp. 183-287, 1877.

\bibitem{LR}P. Lancaster and L. Rodman, ~{\em Algebraic
Riccati Equations},~ Oxford Science Publications, 1995.

%
\bibitem{L1}A. Luger, ``A Factorization of Regular
Generalized Regular Nevanlinna Functions'', {\em Integ.
Eq. \& Op.  Theory}, Vol. 43, pp. 326-345, 2002.
 
\bibitem{Mej}Meijering, ``A Chronology of Interpolation:
From Ancient Astronomy to Modern Signal and Image
Processing", {\em Proc. of the IEEE},~ Vol. 90,
pp. 319-342, 2002.
 
%
%
\bibitem{Ru}W. Rudin,~ {\em Functional Analysis},~
McGraw-Hill, 1973.

\bibitem{Wa}E. Waring, ``Problems Concerning
Interpolations",~ {\em Philos. Tans. Roy. Soc. London},~
Vol. 69, pp. 59-67, 1779.
%
\end{thebibliography}
\end{document}